\documentclass[a4paper,leqno,12pt]{amsart}
\usepackage[leqno]{amsmath}
\usepackage{amstext,amssymb,amsthm, enumitem}
\usepackage{graphicx}
\usepackage{tikz}
\usepackage[cp1250]{inputenc}
\usepackage{float}
\newcommand{\subscript}[2]{$#1 #2$}

\newtheorem*{theorem*}{Theorem}
\newtheorem{theorem}{Theorem}
\newtheorem{proposition}[theorem]{Proposition}
\theoremstyle{definition}
\newtheorem*{prof*}{Proof}
\newtheorem*{proof*}{Proof of Theorem 1}

\tikzset{dots/.append style={ultra thick, fill=none}}

\begin{document}

\title[about weak and restricted weak type inequalities for maximal operators]{On relations between weak and restricted weak type inequalities for maximal operators on non-doubling metric measure spaces}

\author{Dariusz Kosz}
\address{ 
	\newline Faculty of Pure and Applied Mathematics
	\newline Wrocław University of Science and Technology 
	\newline Wyb. Wyspiańskiego 27 
	\newline 50-370 Wrocław, Poland
	\newline \textit{Dariusz.Kosz@pwr.edu.pl}	
}
\begin{abstract} In this article we study a special class of non-doubling metric measure spaces for which there is a significant difference between the incidence of weak and restricted weak type $(p,p)$ inequalities for the centered and non-centered Hardy--Littlewood maximal operators, $M^c$ and $M$. As a corollary we extend the result obtained in \cite{Ko}.
%
%
\end{abstract}

\thanks{
	$\noindent \textit{Key words:}$ Hardy--Littlewood maximal operators, weak and restricted weak type inequalities, non-doubling metric measure spaces.
} 
\maketitle

\section{Introduction}

Consider a metric measure space $\mathbb{X} = (X, \rho, \mu)$ where $\rho$ is a metric and $\mu$ is a Borel measure such that the measure of each ball is finite and strictly positive. By $B(x,r)$ we denote the open ball centered at $x \in X$ with radius $r>0$. If we do not specify the center point and radius we write simply $B$. According to this the $\textit{Hardy--Littlewood}$ $\textit{maximal operators}$, centered $M^c$ and non-centered $M$, are defined respectively by
\begin{displaymath}
M^cf(x) = \sup_{r > 0} \frac{1}{\mu(B(x,r))} \int_{B(x,r)} |f| d\mu, \qquad x \in X,
\end{displaymath}
and
\begin{displaymath}
Mf(x) = \sup_{B \ni x} \frac{1}{\mu(B)} \int_B |f| d \mu , \qquad x \in X.
\end{displaymath} 

Recall that an operator $T$ is said to be of strong type $(p,p)$ for some $p \in [1, \infty]$ if $T$ is bounded on $L^p=L^p(\mathbb{X})$ which means that $ \| T f \|_p \lesssim \| f \|_p$ holds uniformly in $f \in L^p$. Similarly, $T$ is of weak type $(p,p)$ if $T$ is bounded from $L^p$ to $L^{p,\infty}=L^{p,\infty}(\mathbb{X})$ which means that
\begin{equation}
\lambda^p \mu(\{x \colon |Tf(x)|>\lambda\}) \lesssim \| f \|_p
\end{equation}
holds uniformly in $f \in L^p$ and $\lambda >0$ (we use the convention that $L^{\infty,\infty} = L^\infty$). Finally, $T$ is of restricted weak type if $T$ is bounded from $L^{p,1}$ to $L^{p, \infty}$ which in the case $p>1$ is equivalent to the statement that $(1)$ holds uniformly in $f = \chi_E$, $\mu(E) < \infty$, and $\lambda > 0$ (see \cite{BS}, Theorem 5.3, p. 231, for example). It is easy to check that being of strong type $(p,p)$ implies being of weak type $(p,p)$ which in turn implies being of restricted weak type $(p,p)$. Here and anywhere else in this paper the notation $A_1 \lesssim A_2$ is used to indicate that $A_1 \leq CA_2$ with a positive constant $C$ independent of significant quantities.

When dealing with some metric measure space it is usually an important issue to study the mapping properties of the associated maximal operators. First of all, $M^c$ and $M$ are trivially of strong type $(\infty, \infty)$ in case of any metric measure space. It is also a well known fact that if the measure is doubling, that is $\mu(B(x, 2r)) \lesssim \mu(B(x, r))$ uniformly in $x \in X$ and $r > 0$, then they are both of weak type $(1,1)$. However, this statement is not true for metric measure spaces in general.
The next important thing is that the Marcinkiewicz interpolation theorem can be applied for maximal operators in such a way that if $M^c$ (equivalently $M$) is of weak or strong type $(p_0,p_0)$ for some $p_0 \in [1,\infty)$, then it is of strong (and hence weak) type $(p,p)$ for every $p > p_0$. For example, through the interpolation we can deduce that $M^c$ and $M$ are of strong type $(p,p)$ for every $p \in (1, \infty]$ in the case of a doubling space. On the other side there are examples of spaces for which maximal operators are of strong type $(p,p)$ for every $p \in (1, \infty]$ while they are not of weak type $(1,1)$. 
Moreover, there are even examples of spaces for which the associated operators $M^c$ and $M$ are not of weak (and hence strong) type $(p,p)$ for every $p \in [1, \infty)$. All these observations persuade us to describe what we can say about the possible existence of the weak or strong type $(p,p)$ inequalities for $M^c$ and $M$ in general. 

The program of searching spaces with specific mapping properties of the associated maximal operators was greatly contributed by H.-Q. Li. By considering a class of the cusp spaces, he proved several interesting theorems related to this issue (see \cite{Li1}, \cite{Li2} and \cite{Li3}). Lately, the results of H.-Q. Li were extended in \cite{Ko}, where the author characterized all possible configurations of the sets of $p$ for which the weak and strong type $(p,p)$ inequalities for maximal operators, both $M^c$ and $M$, hold. Nevertheless, it is worth noting that all previous works dealing with the mapping properties of maximal operators focused only on weak or strong type $(p,p)$ estimates. The well known fact is that the Marcinkiewicz theorem has a stronger version and to use the interpolation one needs only to show that the maximal operator is of restricted weak type $(p_0, p_0)$ for some $p_0 \in [1,\infty)$ (see \cite{SW}, Theorem 3.15, p. 197, for example). Therefore, a natural step to extend the result obtained in \cite{Ko} is to take into account the restricted weak type inequalities in order to relate them to the weak and strong type inequalities and this is what we do in this article.

\section{Main result}

Let $\mathbb{X}$ be a fixed metric measure space. We denote by $P_s^c$, $P_w^c$ and $P_r^c$ the sets consisting of such $p \in [1, \infty]$ for which the associated operator $M^c$ is of strong, weak, or restricted weak type $(p,p)$, respectively. Similarly, let $P_s$, $P_w$ and $P_r$ consist of such $p \in [1, \infty]$ for which $M$ is of strong, weak, or restricted weak type $(p,p)$, respectively. Then\smallskip

\begin{enumerate}[label=(\roman*)]
	\item each of the six sets is of the form $\{\infty\}$, $[p_0, \infty]$ or $(p_0,\infty]$, for some $p_0 \in [1, \infty)$,\smallskip
	
	\item we have the following inclusions
	\begin{displaymath}
	P_s \subset P_s^c, \quad P_w \subset P_w^c, \quad P_r \subset P_r^c, \quad P_s^c \subset P_w^c \subset  P_r^c \subset \overline{P_s^c}, \quad P_s \subset P_w \subset P_r \subset \overline{P_s},
	\end{displaymath}
	where $\overline{E}$ denotes the closure of $E$ in the usual topology of $\mathbb{R} \cup \{ \infty \}$,
	\item since being of restricted weak and weak type $(1,1)$ is the same, the following implications hold
	\begin{displaymath}
	P_r^c = [1, \infty] \implies P_w^c = [1, \infty], \quad P_r = [1, \infty] \implies P_w = [1, \infty].
	\end{displaymath}
\end{enumerate}

Our intention is to show that $(i)$, $(ii)$ and $(iii)$ are the only necessary conditions that the six sets considered above must satisfy. Namely, we will prove the following.

\begin{theorem}
	Let $P_s^c$, $P_w^c$, $P_r^c$, $P_s$, $P_w$ and $P_r$ be such that the conditions $(i)$, $(ii)$ and $(iii)$ hold. Then there exists a (non-doubling) metric measure space for which the associated Hardy--Littlewood maximal operators, centered $M^c$ and non-centered $M$, satisfy\smallskip
	\begin{itemize}
		\item $M^c$ is of strong type $(p,p)$ if and only if $p \in P_s^c$,\smallskip
		\item $M^c$ is of weak type $(p,p)$ if and only if $p \in P_w^c$,\smallskip
		\item $M^c$ is of restricted weak type $(p,p)$ if and only if $p \in P_r^c$,\smallskip
		\item $M$ is of strong type $(p,p)$ if and only if $p \in P_s$,\smallskip
		\item $M$ is of weak type $(p,p)$ if and only if $p \in P_w$,\smallskip
		\item $M$ is of restricted weak type $(p,p)$ if and only if $p \in P_r$.\smallskip
	\end{itemize} 
\end{theorem}
To prove Theorem 1 we use the technique introduced in \cite{Ko}, where some specific metric measure spaces, namely the first and second generation spaces, were considered. Let us remark several facts relating to these objects:
\begin{enumerate}[label=(\alph*)]
	\item for each space of first generation the equalities $P_s^c = P_s$ and $P_w^c = P_w$ hold,
	\item for each space of second generation $P_s^c = P_w^c = [1, \infty]$, while $P_s$ (and possibly $P_w$) is a proper subset of $[1, \infty]$,
	\item by a suitable mixing of first and second generation spaces we receive a class of spaces characterizing all possible relations between the sets $P_s^c$, $P_w^c$, $P_s$ and $P_w$,
	\item all spaces obtained in such a way in \cite{Ko} satisfy $P_w^c = P_r^c$ and $P_w = P_r$, since the Dirac deltas were used whenever it was shown that the associated $M^c$ or $M$ is not of weak type $(p,p)$ for some $p \in [1, \infty]$.
\end{enumerate}
In this paper we describe other specific spaces, some of which we add to the first generation spaces, and some to the second generation spaces. In particular, we are focused on spaces for which at least one of the equalities $P_w^c = P_r^c$ and $P_w = P_r$ does not occur. For the extended classes of the first and second generation spaces we have the following:

\begin{enumerate}[label=(\subscript{\alph*}{'})]
	\item for each space of first generation $P_s^c = P_s$, $P_w^c = P_w$ and $P_r^c = P_r$ hold,
	\item for each space of second generation $P_s^c = P_w^c = P_r^c = [1, \infty]$, while $P_s$ (and possibly $P_w$ and $P_r$) is a proper subset of $[1, \infty]$,
	\item by a suitable mixing of first and second generation spaces we receive a class of spaces characterizing all possible relations between the sets $P_s^c$, $P_w^c$, $P_r^c$, $P_s$, $P_w$ and $P_r$.
\end{enumerate}

It is worth noting here that the construction of the first and second generation spaces was largely inspired by the work of Stempak, who considered in \cite{St} some specific spaces in the context of modified Hardy--Littlewood maximal operators. The proof of Theorem 1 is located in Section 5.

\section{First generation spaces}
	
	\smallskip First we construct some metric measure spaces (we add them to the class of the first generation spaces defined in \cite{Ko}) for which $P_w^c = P_w = (p_0, \infty]$ while $P_r^c = P_r = [p_0, \infty]$ for any fixed $p_0 \in (1, \infty)$. Note that these conditions imply $P_s^c = P_s = (p_0, \infty]$. We begin with an overview of such spaces and then, after choosing $p_0$, we pass to the details.
	
	Let $\tau = (\tau_{n,i})_{n \in \mathbb{N}, \,  i=1, \dots, n}$ be a fixed system of positive integers satisfying $\frac{\tau_{n,i}}{2^{i-1}} \in \mathbb{N}$. Define 
	\begin{displaymath}
	X_{\tau} = \{x_{n,i,j}, \, x'_{n,i,k} \colon n \in \mathbb{N}, \, i=1, \dots , n, \, j=1, \dots , 2^{i-1}, \, k=1, \dots, \tau_{n,i} \},
	\end{displaymath}
	where all elements $x_{n,i,j}, \, x'_{n,i,k}$ are pairwise different. We use some auxiliary symbols for certain subsets of $X_{\tau}$: for $n \in \mathbb{N}$,
	\begin{displaymath}
	S_n = \{x_{n,i,j}, \, x'_{n,i,k} \colon i=1, \dots , n, \, j=1, \dots , 2^{i-1}, \, k=1, \dots, \tau_{n,i}  \},
	\end{displaymath}
	\begin{displaymath}
	S'_n = \{x'_{n,i,k} \colon i=1, \dots , n, \, k=1, \dots, \tau_{n,i}  \},
	\end{displaymath}
	for $n \in \mathbb{N}$, $i=1, \dots , n$,
	\begin{displaymath}
	S'_{n,i} = \{x'_{n,i,k} \colon k=1, \dots, \tau_{n,i}\},
	\end{displaymath}
	and for $n \in \mathbb{N}$, $1 \leq i 
	\leq i' \leq n$, $j=1, \dots , 2^{i-1}$,  
	\begin{displaymath}
	S'_{n,i',i,j} = \{x'_{n,i',k} \colon k \in (\frac{j-1}{2^{i-1}} \tau_{n,i'}, \frac{j}{2^{i-1}}\tau_{n,i'}]\}.
	\end{displaymath}
	Observe that the sets $S'_{n,i',i,j}$, $j = 1, \dots 2^{i-1}$, are pairwise disjoint, each of them contains exactly $\frac{\tau_{n,i'}}{2^{i-1}}$ elements and $\bigcup_{j=1}^{2^{i-1}} S'_{n,i',i,j} = S'_{n,i'}$. Moreover, for any $1 \leq i_1 \leq i_2 \leq i' \leq n$, $n \in \mathbb{N}$, and $j_l \leq 2^{i_l - 1}$, $l=1,2$, we have either $S'_{n,i',i_2,j_2} \subset S'_{n,i',i_1,j_1}$ or $S'_{n,i',i_1,j_1} \cap S'_{n,i',i_2,j_2} = \emptyset$.
	
	We define the metric $\rho = \rho_\tau$ on $X_\tau$ determining the distance between two different elements $x$ and $y$ by the formula
	\begin{displaymath}
	\rho(x,y) = \left\{ \begin{array}{rl}
	1 & \textrm{if } \{x, y\} = \{x_{n,i,j},x'_{n,i',k}\} \textrm{ and } x'_{n,i',k} \in S'_{n,i',i,j},  \\
	2 & \textrm{in the other case.} \end{array} \right. 
	\end{displaymath}
	Figure 1 shows a model of the space $(X_\tau, \rho)$. The solid line between two points indicates that the distance between them equals $1$. Otherwise the distance equals $2$.  
	
	 \begin{figure}[H]
		\begin{tikzpicture}
		[scale=.8,auto=left,every node/.style={circle,fill,inner sep=2pt}]
		\node[label={[yshift=-1cm]$x_{1,1,1}$}] (n0) at (1,1) {};
		\node[label=$x'_{1,1,1}$] (n1) at (0,4)  {};
		\node[label={[yshift=-0.18cm]$x'_{1,1,\tau_{1,1}}$}] (n2) at (2,4)  {};
		\node[dots] (n3) at (1,4)  {...};
		
		\node[label={[yshift=-1cm]$x_{2,1,1}$}] (c0) at (6,1) {};
		\node[label=$x'_{2,1,1}$] (c1) at (5,4)  {};
		\node[label={[yshift=-0.18cm]$x'_{2,1,\tau_{2,1}}$}] (c2) at (7,4)  {};
		\node[dots] (c3) at (6,4)  {...};
		
		\node[label={[yshift=-1cm]$x_{2,2,1}$}] (l0) at (10,1) {};
		\node[label=$x'_{2,2,1}$] (l1) at (9,4)  {};
		\node[label={[xshift=-0.15cm, yshift=-0.3cm]$x'_{2,2,\tau_{2,2}/2}$}] (l2) at (11,4)  {};
		\node[dots] (l3) at (10,4)  {...};
		
		\node[label={[yshift=-1cm]$x_{2,2,2}$}] (r0) at (14,1) {};
		\node[label={[xshift=0.1cm, yshift=-0.5cm]$x'_{2,2,\tau_{2,2}/2+1}$}] (r1) at (13,4)  {};
		\node[label={[xshift=0.2cm, yshift=-0.2cm]$x'_{2,2,\tau_{2,2}}$}] (r2) at (15,4)  {};
		\node[dots] (r3) at (14,4)  {...};
		
		\node[dots] (oo) at (18,1)  {...};
		
		\foreach \from/\to in {n0/n1, n0/n2, l0/l1, l0/l2, c0/c1, c0/c2, r0/r1, r0/r2, l1/c0, l2/c0, r1/c0, r2/c0}
		\draw (\from) -- (\to);
		\end{tikzpicture}
		\caption{}
	 \end{figure}
	Note that we can explicitly describe any ball: for $n \in \mathbb{N}, \, i=1, \dots , n, \, j=1, \dots , 2^{i-1}$,
	\begin{displaymath}
	B(x_{n,i,j},r) = \left\{ \begin{array}{rl}
	\{x_{n,i,j}\} & \textrm{for } 0 < r \leq 1, \\
	\{x_{n,i,j}\} \cup \bigcup_{i' \geq i} S_{n,i',i,j}& \textrm{for } 1 < r \leq 2,  \\
	X_\tau & \textrm{for } 2 < r, \end{array} \right.
	\end{displaymath} 
	\noindent and for $n \in \mathbb{N}, \, i'=1, \dots , n, \, k=1, \dots, \tau_{n,i}$,
	\begin{displaymath}
	B(x'_{n,i',k},r) = \left\{ \begin{array}{rl}
	\{x'_{n,i',k}\} & \textrm{for } 0 < r \leq 1, \\
	\{x'_{n,i',k}\} \cup \{x_{n, i, j} \colon x'_{n,i',k} \in S_{n,i',i,j}\}& \textrm{for } 1 < r \leq 2,  \\
	X_\tau & \textrm{for } 2 < r. \end{array} \right.
	\end{displaymath} 
	
	We define the measure $\mu = \mu_{\tau, F, m}$ on $X_{\tau}$ by letting $\mu(\{x_{n,i,j}\}) = d_n F(n,i)$ and $\mu(\{x'_{n,i,k}\}) = d_n m_n i$, where $0 <F \leq 1$ is a given function, $m = (m_n)_{n \in \mathbb{N}}$ is a sequence satisfying $m_n \geq 2^n$ and, finally, $d = (d_n)_{n \in \mathbb{N}}$ is a sequence with $d_1 = 1$ and $d_n$ chosen (uniquely for fixed $F$ and $m$) in such a way that $\mu(S_n) = \mu(S_{n-1})/2$, $n \geq 2$ (this implies $\mu(X_\tau) < \infty$). Moreover, observe that $\mu$ is non-doubling. From now on we write simply $|E|$ instead of $\mu(E)$ for $E \subset X_\tau$. The reader should not have any difficulties in identifying when the symbol $|\cdot|$ refers to the measure and when it denotes the absolute value sign.
	
	For a function $ f $ on $ X_\tau $ the Hardy--Littlewood maximal operators, the centered one, $M^c$, and the non-centered one, $M$, are given by
	\begin{displaymath}
	M^cf(x) = \sup_{r > 0} \frac{1}{|B(x,r)|} \sum_{y \in B(x,r)} |f(y)| \cdot |\{y\}|, \qquad x \in X_\tau,
	\end{displaymath}
	and
	\begin{displaymath}
	Mf(x) = \sup_{B \ni x} \frac{1}{|B|} \sum_{y \in B} |f(y)| \cdot |\{y\}|, \qquad x \in X_\tau,
	\end{displaymath}
	respectively. In this setting $M$ is of weak type $(p,p)$ for some $1 < p < \infty$ if $ \| M f \|_{p,\infty} \lesssim \| f \|_p$ uniformly in $f \in \ell^p(X_\tau, \mu)$, where $\| g \|_p = \big( \sum_{x \in X_\tau} |g(x)|^p |\{x\}|\big)^{1/p}$, $ \| g \|_{p,\infty} = \sup_{\lambda > 0} \lambda |E_\lambda(g)|^{1/p}$, and $E_\lambda(g) = \{x \in X_\tau \colon |g(x)| > \lambda\}$. In turn, $M$ is of restricted weak type $(p,p)$ for some $1 < p < \infty$, if $ \| M \chi_E \|_{p,\infty} \lesssim \| \chi_E \|_p$ holds uniformly in $E \subset X_\tau$, where $\chi_E$ is the characteristic function of $E$. Moreover, we introduce the symbol $A_E(f)$ as the average value of a given function $f \geq 0$ on a set $E \subset X_\tau$, namely
	\begin{displaymath}
		A_E(f) = \frac{1}{|E|}\sum_{x \in E} f(x) |\{x\}|. 
	\end{displaymath}
	Analogous definitions and comments apply to $M^c$ instead of $M$ and then to both $M$ and $M^c$ in the context of the space $(Y_\tau, \mu)$ in Section 4. We are ready to describe some of the first generation spaces in detail.
	
	Fix $p_0 \in (1, \infty)$ and consider $\mathbb{X}_{p_0} = (X_\tau, \rho, \mu)$, with a construction based on $m_n = m_n(p_0)$, $\tau_{n,i} = \tau_{n,i}(p_0)$ and $F(n,i) = F_{p_0}(n,i)$, defined as follows. Let $(a_i)_{i \in \mathbb{N}}$ be a sequence satisfying $a_1 = 1$ and $a_i = i^{p_0} - (i-1)^{p_0}$ for $i \geq 2$. For $ i = 1, \dots, n$, $n \in \mathbb{N}$, define
	\begin{displaymath}
	F(n,i)=2^{(i-n)/(p_0-1)}, \qquad \tau_{n,i} = \left \lfloor a_i \right \rfloor 2^{2n\left \lfloor p_0 \right \rfloor} n!/i, \qquad m_n = 2^{(2n\left \lfloor p_0 \right \rfloor -n) / (p_0-1)} (n!)^{1/(p_0-1)}.
	\end{displaymath}
	Observe that $\frac{\tau_{n,i}}{2^{i-1}} \in \mathbb{N}$ and $\tau_{n,i} \, i / m_n^{p_0-1} = 2^n \left \lfloor a_i \right \rfloor$. Moreover, we have $|\{x\}| \geq d_n m_n \geq d_n 2^n \geq |S_n \setminus S_n'|$ for any $x \in S_n'$.

	\begin{proposition}
		Fix $p_0 \in (1, \infty)$ and let $\mathbb{X}_{p_0}$ be the metric measure space defined above. Then the associated maximal operators, centered $M^c$ and non-centered $M$, are not of weak type $(p_0,p_0)$, but are of restricted weak type $(p_0,p_0)$.
	\end{proposition}
	\begin{prof*}
		It suffices to prove that $M^c$ fails to be of weak type $(p_0,p_0)$ and $M$ is of restricted weak type $(p_0,p_0)$. First we show that $M^c$ is not of weak type $(p_0,p_0)$. Indeed, let $f_n = \sum_{i=1}^{n} \sum_{j=1}^{2^{i-1}} 2^{(n-i)/(p_0-1)} \delta_{x_{n,i,j}}$, $n \geq 1$. Then $\|f_n\|_{p_0}^{p_0} = \sum_{i=1}^{n} 2^{i-1} \, 2^{n-i} d_n = 2^{n-1} n d_n$ and
	\begin{displaymath}
	M^cf_n(x'_{n,i',k}) \geq A_{B(x'_{n,i',k}, 3/2)}(f_n) \geq \frac{i' d_n}{2 |\{x'_{n,i',k}\}|} = \frac{1}{2 m_n},
	\end{displaymath}
		for any $x'_{n,i',k} \in S'_n$. This implies that $|E_{1/(4m_n)}(M^cf_n)| \geq |S'_n|$ and hence
	\begin{align*}
\limsup_{n \rightarrow \infty} \frac{\|M^cf_n\|_{p_0,\infty}^{p_0}}{\|f_n\|_{p_0}^{p_0}} &
\geq \lim_{n \rightarrow \infty} \frac{\sum_{i=1}^{n}\tau_{n,i} \, i \, d_n \, m_n}{d_n \, n \, 2^{n-1} \, (4m_n)^{p_0}} 
= 2^{1-2p_0} \lim_{n \rightarrow \infty} \frac{\sum_{i=1}^{n}\tau_{n,i} \, i}{n \, m_n^{p_0-1} \, 2^n} \\
& = 2^{1-2p_0} \lim_{n \rightarrow \infty} \frac{\sum_{i=1}^{n} \left \lfloor a_i \right \rfloor}{n} \geq 2^{-2p_0} \lim_{n \rightarrow \infty} n^{p_0-1} = \infty.
	\end{align*}

In the next step we show that $M$ is of restricted weak type $(p_0,p_0)$. First observe that for any characteristic function $g$ if $\lambda < A_{X_\tau} (g)$, then using Hölder's inequality, we obtain $\lambda^{p_0} |E_\lambda(Mf)| \leq  A_{X_\tau}(g)^{p_0} |X_\tau| \leq \|g\|_{p_0}^{p_0}$. On the other side, with the assumption $\lambda \geq A_{X_\tau} (g)$, if for some $x \in S_n$ we have $Mg(x) > \lambda$, then any ball $B$ containing $x$ and realizing $A_B(g) > \lambda$ must be a subset of $S_n$. This allows us to study the behaviour of $g$ and $Mg$ on each $S_n$ separately. Namely, all that we need is to show that
\begin{equation}
\lambda^{p_0 } | E_\lambda(M\chi_{U})| \lesssim \|\chi_{U}\|_{p_0}^{p_0}
\end{equation}
uniformly in $\lambda > 0$, $\emptyset \neq U \subset S_n$, $n \in \mathbb{N}$.
Fix $n \in \mathbb{N}$ and let $f = \chi_U$ for some $U \subset S_n$. Write $f = \chi_{U'} + \sum_{x \in U \setminus U'} \delta_x$, where $U' = U \cap S_n'$. See that any ball $B_r$ with $B_r \cap S_n \neq \emptyset$ and $r < 2$ is a subset of $S_n$ and for a fixed $i \in \{1, \dots, n\}$ it contains at most one of the points $\{x_{n,i,j} \colon j=1, \dots, 2^{i-1}\}$. This implies that
\begin{displaymath}
Mf(x_0) \leq \max\{A_{X_\tau}(f), \ 2M\chi_{U'}(x_0), \ 2 \, C \max_{x \in U \setminus U'}\{ M\delta_x(x_0)\}\},
\end{displaymath}
with $C = \sum_{i=0}^{\infty} 2^{-i(p_0-1)}$. Therefore it suffices to show that $(2)$ holds uniformly in $\lambda > 0$ (and $n \in \mathbb{N}$) for $U \subset S_n'$ and $U = \{x\} \subset S_n \setminus S_n'$. 

First consider $\lambda > 0$ and $f = \chi_{U}$, where $U \subset S_n'$. Without any loss of generality we can assume that $A_{X_\tau}(f) \leq \lambda < 1$. Therefore $E_\lambda(Mf) \subset S_n$ and since $E_\lambda(Mf) \cap S_n' \neq \emptyset$ we have $|E_\lambda(Mf)| \leq 2 |E_\lambda(Mf) \cap S_n'|$. Consider two different balls, $B_1$ and $B_2$, both contained in $S_n$, and denote $B_l' = B_l \cap S_n'$, $l=1,2$. Then we have one of the three possibilities: $B_1' \subset B_2'$, $B_2' \subset B_1'$ or $B_1' \cap B_2' = \emptyset$. Combining this observation with the fact that $f(x) = 0$ for $x \notin S_n'$, we obtain $A_{E_\lambda(Mf) \cap S_n'}(f) > \lambda$ and hence
\begin{displaymath}
\lambda^{p_0} |E_\lambda(Mf)| < 2 A_{E_\lambda(Mf) \cap S_n'}(f)^{p_0} \, |E_\lambda(Mf) \cap S_n'| \leq 2 \|f\|_{p_0}^{p_0}.
\end{displaymath}

Finally, consider $\lambda > 0$ and $\delta_{x_{n,i,j}}$ for some fixed $x_{n,i,j} \in S_n \setminus S_n'$. We may assume that $A_{X_\tau}(\delta_{x_{n,i,j}}) \leq \lambda < 1$. If $E_\lambda(M\delta_{x_{n,i,j}}) \cap S_n' = \emptyset$ then $E_\lambda(M\delta_{x_{n,i,j}}) = \{x_{n,i,j}\}$ and $\lambda^{p_0} |E_\lambda(M\delta_{x_{n,i,j}})| \leq \|\delta_{x_{n,i,j}}\|_{p_0}^{p_0}$ holds trivially. Otherwise, if $E_\lambda(M\delta_{x_{n,i,j}}) \cap S_n' \neq \emptyset$, then $|E_\lambda(M\delta_{x_{n,i,j}})| \leq 2 |E_\lambda(M\delta_{x_{n,i,j}}) \cap S_n'|$. For any $x'_{n,i',k} \in S'_{n,i',i,j}$ we have the estimate
\begin{displaymath}
M\delta_{x_{n,i,j}}(x'_{n,i',k}) = A_{B(x'_{n,i',k}, 3/2)}(\delta_{x_{n,i,j}}) \leq \frac{2^{(i-n)/(p_0-1)}}{ m_n \, i'},
\end{displaymath}
while $M\delta_{x_{n,i,j}}(x) = A_{X_\tau}(\delta_{x_{n,i,j}})$ for $x \in S_n' \setminus \bigcup_{i'=i}^{n} S'_{n,i',i,j}$. For any $i' \in \{i, \dots, n\}$ we obtain
\begin{align*}
\Big(   \frac{2^{(i-n)/(p_0-1)}}{ m_n \, i'}   \Big)^{p_0} \sum_{l=i}^{i'} |S'_{n,l,i,j}| &= 
\frac{2^{(i-n)p_0 / (p_0-1)}}{( m_n \, i')^{p_0}}    \sum_{l=i}^{i'} \frac{|S'_{n,l}|}{2^{i-1}} = 
\frac{2^{(i-n)p_0 / (p_0-1)}}{i'^{p_0}}  \sum_{l=i}^{i'} \frac{d_n \, \tau_{n,l} \, l}{m_n^{p_0-1} \, 2^{i-1}} \\
& = \frac{d_n \, 2^{(i-n)p_0 / (p_0-1)}}{i'^{p_0}} \sum_{l=i}^{i'} 2^{n-i+1} \left \lfloor a_l \right \rfloor \\ & \leq d_n \, 2^{((i-n)/(p_0-1)) + 1} \, \frac{\sum_{l=1}^{i'} a_l}{i'^{p_0}} = 2 \|\delta_{x_{n,i,j}}\|_{p_0}^{p_0},
\end{align*}
and hence $\lambda^{p_0} |E_\lambda(M\delta_{x_{n,i,j}})| \leq 4 \|\delta_{x_{n,i,j}}\|_{p_0}^{p_0}$ uniformly in $\lambda > 0$ and $x_{n,i,j} \in X_\tau$. $\raggedright \hfill \qed$
		
\end{prof*}

\section{Second generation spaces}

Now we construct some metric measure spaces (we add them to the class of the second generation spaces defined in \cite{Ko}) for which $P_w^c = P_r^c = [1, \infty)$, $P_w = (p_0, \infty]$ and $P_r = [p_0, \infty]$ with any fixed $p_0 \in (1, \infty)$. Let $\tau = (\tau_{n,i})_{n \in \mathbb{N}, \,  i=1, \dots, n}$ be a fixed system of positive integers satisfying $\frac{\tau_{n,i}}{2^{i-1}} \in \mathbb{N}$. Define 
\begin{displaymath}
Y_{\tau} = \{y_{n,i,j}, \, y^\circ_{n,i,k}, \, y'_{n,i,k} \colon n \in \mathbb{N}, \, i=1, \dots , n, \, j=1, \dots , 2^{i-1}, \, k=1, \dots, \tau_{n,i} \},
\end{displaymath}
where all elements $y_{n,i,j}, \, y^\circ_{n,i,k}, \, y'_{n,i,k}$ are pairwise different. We use some auxiliary symbols for certain subsets of $X_{\tau}$: for $n \in \mathbb{N}$,
\begin{displaymath}
T_n = \{y_{n,i,j}, \, y^\circ_{n,i,k}, \, y'_{n,i,k} \colon i=1, \dots , n, \, j=1, \dots , 2^{i-1}, \, k=1, \dots, \tau_{n,i}  \},
\end{displaymath}
\begin{displaymath}
T^\circ_n = \{y^\circ_{n,i,k} \colon i=1, \dots , n, \, k=1, \dots, \tau_{n,i}  \}, \quad T'_n = \{y'_{n,i,k} \colon i=1, \dots , n, \, k=1, \dots, \tau_{n,i}  \},
\end{displaymath}
for $n \in \mathbb{N}$, $i=1, \dots , n$,
\begin{displaymath}
T^\circ_{n,i} = \{y^\circ_{n,i,k} \colon k=1, \dots, \tau_{n,i}\}, \quad T'_{n,i} = \{y'_{n,i,k} \colon k=1, \dots, \tau_{n,i}\},
\end{displaymath}
and for $j=1, \dots , 2^{i-1}$, $1 \leq i 
\leq i' \leq n$, $n \in \mathbb{N}$,
\begin{displaymath}
T^\circ_{n,i',i,j} = \{y^\circ_{n,i',k} \colon k \in (\frac{j-1}{2^{i-1}} \tau_{n,i'}, \frac{j}{2^{i-1}}\tau_{n,i'}]\}, \quad
T'_{n,i',i,j} = \{y'_{n,i',k} \colon k \in (\frac{j-1}{2^{i-1}} \tau_{n,i'}, \frac{j}{2^{i-1}}\tau_{n,i'}]\}.
\end{displaymath}

We define the metric $\rho = \rho_\tau$ on $Y_\tau$ determining the distance between two different elements $x$ and $y$ by the formula
\begin{displaymath}
\rho(x,y) = \left\{ \begin{array}{rl}
1 & \textrm{if } \{x, y\} = \{y_{n,i,j},y^\circ_{n,i',k}\} \textrm{ and } x^\circ_{n,i',k} \in T^\circ_{n,i',i,j}, \\
1 & \textrm{if } \{x, y\} \subset T_n \setminus (T_n^\circ \cup T_n'), \\
1 & \textrm{if } \{x, y\} = \{y^\circ_{n, i', k}, y'_{n,i',k}\}, \\
2 & \textrm{in the other case.} \end{array} \right. 
\end{displaymath}

	Figure 2 shows a model of the space $(Y_\tau, \rho)$ (with the convention as in Figure 1).
	
	\begin{figure}[H]
		\begin{tikzpicture}
		[scale=.8,auto=left,every node/.style={circle,fill,inner sep=2pt}]
		
		\node[label={[yshift=-1cm]$y_{1,1,1}$}] (n0) at (1,1) {};
		
		\node[label=$y^\circ_{1,1,1}$] (n1) at (0,4)  {};
		\node[label={[yshift=-0.18cm]$y^\circ_{1,1,\tau_{1,1}}$}] (n2) at (2,4)  {};
		\node[dots] (n3) at (1,4)  {...};
		
		\node[label=$y'_{1,1,1}$] (n4) at (0,7)  {};
		\node[label={[yshift=-0.18cm]$y'_{1,1,\tau_{1,1}}$}] (n5) at (2,7)  {};
		\node[dots] (n6) at (1,7)  {...};

		\node[label={[yshift=-1cm]$y_{2,1,1}$}] (c0) at (6,1) {};
		
		\node[label=$y^\circ_{2,1,1}$] (c1) at (5,4)  {};
		\node[label={[yshift=-0.18cm]$y^\circ_{2,1,\tau_{2,1}}$}] (c2) at (7,4)  {};
		\node[dots] (c3) at (6,4)  {...};
		
		\node[label=$y'_{2,1,1}$] (c4) at (5,7)  {};
		\node[label={[yshift=-0.18cm]$y'_{2,1,\tau_{2,1}}$}] (c5) at (7,7)  {};
		\node[dots] (c6) at (6,7)  {...};

		\node[label={[yshift=-1cm]$y_{2,2,1}$}] (l0) at (10,1) {};
		
		\node[label=$y^\circ_{2,2,1}$] (l1) at (9,4)  {};
		\node[label={[xshift=-0.15cm, yshift=-0.3cm]$y^\circ_{2,2,\tau_{2,2}/2}$}] (l2) at (11,4)  {};
		\node[dots] (l3) at (10,4)  {...};
		
		\node[label=$y'_{2,2,1}$] (l4) at (9,7)  {};
		\node[label={[xshift=-0.15cm, yshift=-0.3cm]$y'_{2,2,\tau_{2,2}/2}$}] (l5) at (11,7)  {};
		\node[dots] (l6) at (10,7)  {...};

		\node[label={[yshift=-1cm]$y_{2,2,2}$}] (r0) at (14,1) {};
		
		\node[label={[xshift=0.1cm, yshift=-0.5cm]$y^\circ_{2,2,\tau_{2,2}/2+1}$}] (r1) at (13,4)  {};
		\node[label={[xshift=0.2cm, yshift=-0.2cm]$y^\circ_{2,2,\tau_{2,2}}$}] (r2) at (15,4)  {};
		\node[dots] (r3) at (14,4)  {...};
		
		\node[label={[xshift=0.1cm, yshift=-0.5cm]$y'_{2,2,\tau_{2,2}/2+1}$}] (r4) at (13,7)  {};
		\node[label={[xshift=0.2cm, yshift=-0.2cm]$y'_{2,2,\tau_{2,2}}$}] (r5) at (15,7)  {};
		\node[dots] (r6) at (14,7)  {...};

		\node[dots] (oo) at (18,1)  {...};

		\foreach \from/\to in {n0/n1, n0/n2, l0/l1, l0/l2, c0/c1, c0/c2, r0/r1, r0/r2, l1/c0, l2/c0, r1/c0, r2/c0, l0/c0, l0/r0}
		\draw (\from) -- (\to);
		
		\draw (0,7) -- (0,5); \draw (0,4.5) -- (0,4);
		\draw (2,7) -- (2,5); \draw (2,4.5) -- (2,4);
		\draw (5,7) -- (5,5); \draw (5,4.5) -- (5,4);
		\draw (7,7) -- (7,5); \draw (7,4.5) -- (7,4);
		\draw (9,7) -- (9,5); \draw (9,4.5) -- (9,4);
		\draw (11,7) -- (11,5); \draw (11,4.5) -- (11,4);
		\draw (13,7) -- (13,5); \draw (13,4.5) -- (13,4);
		\draw (15,7) -- (15,5); \draw (15,4.5) -- (15,4);
		\draw (6,1) arc (-113.75:-66.25:10);
		\end{tikzpicture}
		\caption{}
	\end{figure}
	Note that we can explicitly describe any ball: for $n \in \mathbb{N}, \, i=1, \dots , n, \, j=1, \dots , 2^{i-1}$,
	\begin{displaymath}
	B(y_{n,i,j},r) = \left\{ \begin{array}{rl}
	\{y_{n,i,j}\} & \textrm{for } 0 < r \leq 1, \\
	\big( T_n \setminus (T_n^\circ \cup T_n') \big) \, \cup \, \bigcup_{i' \geq i} T^\circ_{n,i',i,j}& \textrm{for } 1 < r \leq 2,  \\
	Y_\tau & \textrm{for } 2 < r, \end{array} \right.
	\end{displaymath} 
	\noindent and for $n \in \mathbb{N}, \, i'=1, \dots , n, \, k=1, \dots, \tau_{n,i}$,
	\begin{displaymath}
	B(y^\circ_{n,i',k},r) = \left\{ \begin{array}{rl}
	\{y^\circ_{n,i',k}\} & \textrm{for } 0 < r \leq 1, \\
	\{y^\circ_{n,i',k}, y'_{n,i',k}\} \cup \{y_{n, i, j} \colon y^\circ_{n,i',k} \in T^\circ_{n,i',i,j}\}& \textrm{for } 1 < r \leq 2,  \\
	Y_\tau & \textrm{for } 2 < r, \end{array} \right.
	\end{displaymath} 
	\noindent and 
	\begin{displaymath}
	B(y'_{n,i',k},r) = \left\{ \begin{array}{rl}
	\{y'_{n,i',k}\} & \textrm{for } 0 < r \leq 1, \\
	\{y^\circ_{n,i',k}, y'_{n,i',k}\} & \textrm{for } 1 < r \leq 2,  \\
	Y_\tau & \textrm{for } 2 < r. \end{array} \right.
	\end{displaymath} 
	
	We define the measure $\mu = \mu_{\tau, F, G, m}$ on $Y_{\tau}$ by letting $\mu(\{y_{n,i,j}\}) = d_n F(n,i)$, $\mu(\{y^\circ_{n,i,j}\}) = d_n G(n)$, $\mu(\{x'_{n,i,j}\}) = d_n m_n i$, where $0 <F \leq 1$ and $0 < G \leq 1 / \sum_{i=1}^{n} \tau_{n,i}$ are given functions, $m = (m_n)_{n \in \mathbb{N}}$ is a sequence satisfying $m_n \geq 2^n$ and finally $d = (d_n)_{n \in \mathbb{N}}$ is a sequence with $d_1 = 1$ and $d_n$ chosen (uniquely for fixed $F$, $G$ and $m$) in such a way that $\mu(T_n) = \mu(T_{n-1})/2$, $n \geq 2$.
	
	Now let $p_0 \in (1, \infty)$, and consider $\mathbb{Y}_{p_0} = (Y_\tau, \rho, \mu)$, with a construction based on $m$, $\tau$, $F$ defined as in Section 3 and $G(n) = 2^{(1-n)/(p_0-1)} / \sum_{i=1}^{n}\tau_{n,i}$. Observe that we have $|\{y\}| \geq d_n m_n \geq d_n 2^n \geq |T_n \setminus T_n'|$ for any $y \in T_n'$. Moreover, we have $|\{y\}| \geq |T_n^\circ|$ for any $y \in T_n \setminus (T_n^\circ \cup T_n')$.
	
	\begin{proposition}
		Fix $p_0 \in (1, \infty)$ and let $\mathbb{Y}_{p_0}$ be the metric measure space defined above. Then the associated centered maximal operator $M^c$ is of strong type $(1,1)$ while the non-centered $M$ is not of weak type $(p_0,p_0)$, but is of restricted weak type $(p_0,p_0)$.
	\end{proposition}
	
	\begin{prof*}
		First we show that $M^c$ is of strong type $(1,1)$. Let $f \in \ell^1(\hat{\mathbb{Y}}_{p_0})$, $f \geq 0$. We use the estimate: for $n \in \mathbb{N}, \, i=1, \dots , n, \, j=1, \dots , 2^{i-1}$,
		\begin{displaymath}
		Mf(y_{n,i,j}) \leq f(y_{n,i,j}) + 2 A_{T_n \setminus T_n'}(f) + A_{Y_\tau}(f),
		\end{displaymath}
		\noindent and for $n \in \mathbb{N}, \, i'=1, \dots , n, \, k=1, \dots, \tau_{n,i}$,
		\begin{displaymath}
		Mf(y^\circ_{n,i',k}) \leq f(y^\circ_{n,i',k}) + \sup_{y \in T_n \setminus T_n^\circ}f(y) + A_{Y_\tau}(f),
		\end{displaymath}
		\noindent and
		\begin{displaymath}
		Mf(y'_{n,i',k}) \leq f(y'_{n,i',k}) + A_{\{y^\circ_{n, i', k}, y'_{n, i',k}\}}(f) + A_{Y_\tau}(f).
		\end{displaymath}
		Observe that 
		\begin{displaymath}
		\sum_{n \in \mathbb{N}}2 A_{T_n \setminus T_n'}(f) \cdot |T_n \setminus (T_n^\circ \cup T_n')| \leq 2 \|f\|_1,
		\end{displaymath} 
		and 
		\begin{displaymath}
		\sum_{n \in \mathbb{N}} \sum_{i'=1}^n \sum_{k=1}^{\tau_{n,i}} A_{\{y^\circ_{n, i', k}, y'_{n, i',k}\}} \cdot |\{y'_{n, i',k}\}| \leq \|f\|_1.
		\end{displaymath} 
		Moreover, since $|\{y\}| \geq |T_n^\circ|$ for any $y \in T_n \setminus T_n^\circ$, we have
		\begin{displaymath}
		\sum_{n \in \mathbb{N}} \sum_{i'=1}^n \sum_{k=1}^{\tau_{n,i}} \sup_{y \in T_n \setminus T_n^\circ}f(y) \cdot |\{y^\circ_{n,i',k}\}| \leq \sum_{n \in \mathbb{N}} \sup_{y \in T_n \setminus T_n^\circ} \Big( f(y) \cdot |\{y\}| \Big) \leq \|f\|_1,
		\end{displaymath}
		and hence we obtain $\|Mf\|_1 \leq 6 \|f\|_1$.
		
		In the next step we show that $M$ is not of weak type $(p_0, p_0)$. Indeed, let $f_n = \sum_{i=1}^{n} \sum_{j=1}^{2^{i-1}} 2^{(n-i)/(p_0-1)} \delta_{y_{n,i,j}}$, $n \geq 1$. Then $\|f_n\|_{p_0}^{p_0} = 2^{n-1} n d_n$ and
		\begin{displaymath}
		Mf_n(y'_{n,i',k}) \geq A_{B(y^\circ_{n,i',k}, 3/2)}(f_n) \geq \frac{i' d_n}{2 |\{y'_{n,i',k}\}|} = \frac{1}{2 m_n},
		\end{displaymath}
		for any $y'_{n,i',k} \in T'_n$. This implies that $|E_{1/(4m_n)}(Mf_n)| \geq |T'_n|$ and hence
		\begin{align*}
		\limsup_{n \rightarrow \infty} \frac{\|Mf_n\|_{p_0,\infty}^{p_0}}{\|f_n\|_{p_0}^{p_0}} &
		\geq \lim_{n \rightarrow \infty} \frac{\sum_{i=1}^{n}\tau_{n,i} \, i \, d_n \, m_n}{d_n \, n \, 2^{n-1} \, (4m_n)^{p_0}} 
		= 2^{1-2p_0} \lim_{n \rightarrow \infty} \frac{\sum_{i=1}^{n}\tau_{n,i} \, i}{n \, m_n^{p_0-1} \, 2^n} \\
		& = 2^{1-2p_0} \lim_{n \rightarrow \infty} \frac{\sum_{i=1}^{n} \left \lfloor a_i \right \rfloor}{n} \geq 2^{-2p_0} \lim_{n \rightarrow \infty} n^{p_0-1} = \infty.
		\end{align*}
		
		In the last step we show that $M$ is of restricted weak type $(p_0,p_0)$. Arguing similarly as in the proof of Proposition 1 we observe that it suffices to show that
		\begin{equation}
		\lambda^{p_0 } | E_\lambda(M\chi_{U})| \lesssim \|\chi_{U}\|_{p_0}^{p_0}
		\end{equation}
		uniformly in $\lambda > 0$, $\emptyset \neq U \subset T_n$, $n \in \mathbb{N}$.
		Fix $n \in \mathbb{N}$ and let $f = \chi_U$ for some $U \subset T_n$. Write $f = \chi_{U^\circ \cup U'} + \sum_{y \in U \setminus (U^\circ \cup U')} \delta_y$, where $U^\circ = U \cap T_n^\circ$ and $U' = U \cap T_n'$. For any $y_0 \in T_n$ we have the estimate
		\begin{displaymath}
		Mf(y_0) \leq \max\{A_{Y_\tau}(f), \ 3M\chi_{U \setminus (U^\circ \cup U')}(y_0), \ 3M\chi_{U^\circ}(y_0), \ 3M\chi_{U'}(y_0)\}, 
		\end{displaymath}
		and hence it suffices to show that $(3)$ holds uniformly in $\lambda > 0$ (and $n \in \mathbb{N}$) for $U$ being a subset of $T_n \setminus (T_n^\circ \cup T_n')$, $T_n^\circ$ or $T_n'$. Moreover, note that any ball $B_r$ with $B_r \cap T_n' \neq \emptyset$ and $r < 2$ is a subset of $T_n$ and for a fixed $i \in \{1, \dots, n\}$ it contains at most one of the points $\{y_{n,i,j} \colon j=1, \dots, 2^{i-1}\}$. This implies that for $y_0 \in T_n'$ we have
		\begin{displaymath}
		Mf(y_0) \leq \max\{A_{Y_\tau}(f), \ 3M\chi_{U^\circ}(y_0), \ 3M\chi_{U'}(y_0), \ 3 \, C \max_{y \in U \setminus (U^\circ \cup U')}\{ M\delta_y(y_0)\}\},
		\end{displaymath}
		with $C = \sum_{i=0}^{\infty} 2^{-i(p_0-1)}$. Recall that $|\{y\}| \geq d_n m_n \geq d_n 2^n \geq |T_n \setminus T_n'|$ for any $y \in T_n'$. Combining these two observations we conclude that if $E_\lambda(Mf) \cap T_n' \neq \emptyset$, then it suffices to show that $(3)$ holds uniformly in $\lambda > 0$ (and $n \in \mathbb{N}$) for $U \subset T_n'$, $U \subset T_n^\circ$, or $U = \{y\} \subset T_n \setminus (T_n^\circ \cup T_n')$.
		
		With this conclusion consider first $f = \chi_{U}$ such that $U \subset T_n \setminus (T_n^\circ \cup T_n')$ and $\lambda > 0$ such that $A_{Y_\tau}(f) \leq \lambda < 1$. We have one of the two possibilities: $E_\lambda(Mf) \cap T_n' = \emptyset$ or $E_\lambda(Mf) \cap T_n' \neq \emptyset$. If $E_\lambda(Mf) \cap T_n' \neq \emptyset$, then, without any loss of generality, we can take $U = \{y_{n,i,j}\}$. For any $y'_{n,i',k} \in T'_{n,i',i,j}$ we have the estimate
		\begin{displaymath}
		M\delta_{y_{n,i,j}}(y'_{n,i',k}) = A_{B(y^\circ_{n,i',k}, 3/2)}(\delta_{y_{n,i,j}}) \leq \frac{2^{(i-n)/(p_0-1)}}{ m_n \, i'},
		\end{displaymath}
		while $M\delta_{y_{n,i,j}}(y) = A_{Y_\tau}(\delta_{y_{n,i,j}})$ for $y \in T_n' \setminus \bigcup_{i'=i}^{n} T'_{n,i',i,j}$. For any $i' \in \{i, \dots, n\}$ we obtain
		\begin{align*}
		\Big(   \frac{2^{(i-n)/(p_0-1)}}{ m_n \, i'}   \Big)^{p_0} \sum_{l=i}^{i'} |T'_{n,l,i,j}| &= 
		\frac{2^{(i-n)p_0 / (p_0-1)}}{( m_n \, i')^{p_0}}    \sum_{l=i}^{i'} \frac{|T'_{n,l}|}{2^{i-1}} = 
		\frac{2^{(i-n)p_0 / (p_0-1)}}{i'^{p_0}}  \sum_{l=i}^{i'} \frac{d_n \, \tau_{n,l} \, l}{m_n^{p_0-1} \, 2^{i-1}} \\
		& = \frac{d_n \, 2^{(i-n)p_0 / (p_0-1)}}{i'^{p_0}} \sum_{l=i}^{i'} 2^{n-i+1} \left \lfloor a_l \right \rfloor \\ &\leq d_n \, 2^{((i-n)/(p_0-1)) + 1} \, \frac{\sum_{l=1}^{i'} a_l}{i'^{p_0}} = 2 \|\delta_{y_{n,i,j}}\|_{p_0}^{p_0},
		\end{align*}
		and hence $\lambda^{p_0} |E_\lambda(M\delta_{y_{n,i,j}})| \leq 4 \|\delta_{y_{n,i,j}}\|_{p_0}^{p_0}$. 
		Therefore, consider the case $E_\lambda(Mf) \cap T_n' = \emptyset$. Since $E_\lambda(Mf) \cap T_n \setminus (T_n^\circ \cup T_n') \neq \emptyset$, we have $|E_\lambda(Mf)| \leq 2 |E_\lambda(Mf) \cap T_n \setminus (T_n^\circ \cup T_n')|$. If $\lambda \leq A_{T_n \setminus (T_n^\circ \cup T_n')}(f)$, then
		\begin{displaymath}
			\lambda^{p_0 } | E_\lambda(Mf)| \leq 2 \lambda^{p_0 } |E_\lambda(Mf) \cap T_n \setminus (T_n^\circ \cup T_n')| \leq 2 \|f\|_{p_0}^{p_0}.
		\end{displaymath}
		Otherwise, let $\lambda > A_{T_n \setminus (T_n^\circ \cup T_n')}(f)$ and fix $y \in E_\lambda(Mf) \cap T_n \setminus (T_n^\circ \cup T_n')$. Observe that the volume of each ball $B_r$ such that $y \in B_r$ and $r > 1$, is greater than $|T_n \setminus (T_n^\circ \cup T_n')|$ and hence $A_{B_r}(f) \leq \lambda$. This implies that $y \in U$ and therefore 
		\begin{displaymath}
			\lambda^{p_0 } | E_\lambda(Mf)| \leq 2 |E_\lambda(Mf) \cap T_n \setminus (T_n^\circ \cup T_n')| = 2 \|f\|_{p_0}^{p_0}.
		\end{displaymath}
		
		Next, consider $f = \chi_{U}$ such that $U \subset T_n'$ and $\lambda > 0$ such that $A_{Y_\tau}(f) \leq \lambda < 1$. Then $E_\lambda(Mf) \subset T_n$ and, since $E_\lambda(Mf) \cap T_n' \neq \emptyset$, we have $|E_\lambda(Mf)| \leq 2 |E_\lambda(Mf) \cap T_n'|$. Observe that there is no ball $B \subset T_n$ which contains two different points from $T_n'$. Therefore, if $y \in E_\lambda(Mf) \cap T_n'$, then $y \in U$ and hence $\lambda^{p_0 } | E_\lambda(Mf)| \leq  2\lambda^{p_0 } |E_\lambda(Mf) \cap T_n'| \leq 2 \|f\|_{p_0}^{p_0}$.
		
		Lastly, consider $f = \chi_{U}$ such that $U \subset T_n^\circ$ and $\lambda > 0$ such that $A_{Y_\tau}(f) \leq \lambda < 1$. Assume that $E_\lambda(Mf) \cap T_n' \neq \emptyset$. If $y'_{n,i',k} \in E_\lambda(Mf) \cap T_n'$, then $y^\circ_{n,i',k} \in U$ and $\lambda < Mf(y'_{n,i',k}) = A_{B(y'_{n,i',k}, 3/2)}(f) \leq |\{y^\circ_{n,i',k}\}| / |\{y'_{n,i',k}\}|$. Therefore we obtain
		\begin{displaymath}
		\lambda^{p_0 } | E_\lambda(Mf)| \leq 2\lambda^{p_0 } |E_\lambda(Mf)\cap T_n'| \leq 2 \|f\|_{p_0}^{p_0}.
		\end{displaymath}
		On the other hand, assume that $E_\lambda(Mf) \subset T_n \setminus T_n'$ and  $E_\lambda(Mf) \not\subset T_n^\circ$. Hence we have $|E_\lambda(Mf)| \leq 2 |E_\lambda(Mf) \cap T_n \setminus (T_n^\circ \cup T_n')|$. Moreover, the volume of any ball containing points $x$ and $y$ such that $x \in T_n \setminus (T_n^\circ \cup T_n')$ and $y \in T_n^\circ$, is greater than $|T_n \setminus (T_n^\circ \cup T_n')|$. Therefore, if $y \in E_\lambda(Mf) \cap T_n \setminus (T_n^\circ \cup T_n')$, then $\lambda < Mf(y) < \|f\|_1 / |T_n \setminus (T_n^\circ \cup T_n')|$, and using the similar argument as before we obtain $\lambda^{p_0 } | E_\lambda(Mf)| \leq 2 \|f\|_{p_0}^{p_0}.$
		Finally, assume the last case which is $E_\lambda(Mf) \subset T_n^\circ$. Because there is no ball $B \subset T_n^\circ$ which contains two different points from $T_n^\circ$, then $E_\lambda(Mf) = U$ and therefore $\lambda^{p_0} |E_\lambda(Mf)| \leq \|f\|_{p_0}^{p_0}. \raggedright \hfill \qed$
	\end{prof*}
	
	\section{Proof of Theorem 1}
	
	The specific technique used to construct all spaces discussed earlier in this paper and in \cite{Ko} ensures that the possible interactions between the different branches are completely meaningless while studying the existence of the restricted weak, weak and strong type inequalities. This fact plays a crucial role in the proof of Theorem 1, where we want to obtain a wider spectrum of possible behaviours of $M^c$ and $M$ by a suitable mixing of the first and second generation spaces. We will explain this idea in detail shortly.
	
	\begin{proof*}
		We consider only two cases. If the equalities $P_w^c = P_r^c$ and $P_w = P_r$ hold, then the expected space may be chosen in such a way as it was done in \cite{Ko} for the sets $P_s^c$, $P_w^c$, $P_s$, and $P_w$. Assume that the opposite is true. Note that in this case $P_s$ must be a proper subset of $[1, \infty]$. We can find a space $\mathbb{X} = (X, \rho_X, \mu_X)$ of first generation for which 
		\begin{displaymath}
		P_s^c(\mathbb{X}) = P_s(\mathbb{X}) = P_s^c, \quad P_w^c(\mathbb{X}) = P_w(\mathbb{X}) = P_w^c, \quad
		P_r^c(\mathbb{X}) = P_r(\mathbb{X}) = P_r^c,
		\end{displaymath}
		hold (if $P_w^c = P_r^c$, then the desired space is one of the spaces described in Section 2 in \cite{Ko}; otherwise, if $P_w^c \neq P_r^c$, then the conditions $(i)-(iii)$ imply that $P_w^c = (p_0, \infty]$ and $P_r^c = [p_0, \infty]$ for some $p_0 \in (1, \infty)$ and hence the desired space is one of the spaces described in Section 3 in this paper). Similarly, regardless of the possibilities, $P_w = P_r$ or $P_w \neq P_r$, we can find a space $\mathbb{Y} = (Y, \rho_Y, \mu_Y)$ of second generation for which
		\begin{displaymath}
		P_s^c(\mathbb{Y}) = P_w^c(\mathbb{Y}) = P_r^c(\mathbb{Y}) = [1, \infty], \quad
		P_s(\mathbb{Y}) = P_s, \quad P_w(\mathbb{Y}) = P_w, \quad P_r(\mathbb{Y}) = P_r.
		\end{displaymath}
		
		Using $\mathbb{X}$ and $\mathbb{Y}$ and assuming that $X \cap Y = \emptyset$ we construct the space $\mathbb{Z} = (Z, \rho_Z, \mu_Z)$ as follows.
		Denote $Z = X \cup Y$. We define the metric $\rho_Z$ on $Z$ by
		\begin{displaymath}
		\rho_Z(x,y) = \left\{ \begin{array}{rl}
		\rho_X(x,y) & \textrm{if }  \{x,y\} \subset X,   \\
		\rho_Y(x,y) & \textrm{if }  \{x,y\} \subset Y,   \\
		2 & \textrm{in the other case,} \end{array} \right. 
		\end{displaymath} 
		and the measure $\mu_Z$ on $Z$ by
		\begin{displaymath}
		\mu_Z(E) = \mu_X(E \cap X) + \mu_Y(E \cap Y), \qquad E \subset Z.
		\end{displaymath}
		It can easily be shown that $\mathbb{Z}$ has the following properties
		\begin{itemize}
			\item $P_s^c(\mathbb{Z})=P_s^c(\mathbb{X}) \cap P_s^c(\mathbb{Y})= P_s^c \cap [1, \infty] = P_s^c$, \smallskip
			\item $P_w^c(\mathbb{Z})=P_w^c(\mathbb{X}) \cap P_w^c(\mathbb{Y})= P_w^c \cap [1, \infty] = P_w^c$, \smallskip
			\item $P_r^c(\mathbb{Z})=P_r^c(\mathbb{X}) \cap P_r^c(\mathbb{Y})= P_r^c \cap [1, \infty] = P_r^c$, \smallskip
			\item $P_s(\mathbb{Z}) = P_s(\mathbb{X}) \cap P_s(\mathbb{Y}) = P_s^c \cap P_s = P_s$, \smallskip
			\item $P_w(\mathbb{Z}) = P_w(\mathbb{X}) \cap P_w(\mathbb{Y}) = P_w^c \cap P_w = P_w$, \smallskip
			\item $P_r(\mathbb{Z}) = P_r(\mathbb{X}) \cap P_r(\mathbb{Y}) = P_r^c \cap P_r = P_r$, 
		\end{itemize} 
		and therefore $\mathbb{Z}$ may be chosen to be the expected space.
		Finally, it is not hard to see that $\mu_Z$ is non-doubling, since it is bounded and there is a ball $B$ in $Z$ with radius $r = 1$ and $|B| < \epsilon$ for any arbitrarily small $\epsilon > 0$.
		$\raggedright \hfill \qed$
	\end{proof*}
	
	\section*{Acknowledgement}
	I would like to thank my supervisor Professor Krzysztof Stempak who posed the problem discussed in this article. His attitude was a great motivation and help for me to achieve the pointed target. I am also indebted to the referee for a very careful reading of the manuscript and suggesting many valuable improvements.
	
	Research was supported by the National Science Centre of Poland, project no. \linebreak
	2016/21/N/ST1/01496.

\end{document}